\documentclass[12pt]{article}
 \usepackage{amssymb}
   \usepackage{graphics}
    \usepackage[dvips]{epsfig}
   \usepackage{times}
   \usepackage{tikz}
   \usepackage{tcolorbox}
   \usepackage{mathtools}
   \newcommand{\D}{\displaystyle}
   
    \def\B{{\mathbb B}}
   \def\C{{\mathbb C}}
   
   \def\H{{\mathbb H}}

   \def\R{{\mathbb R}}
   \def\N{{\mathbb N}}

   \def\Pf{{\it Proof.$\;\;$}}

   \def\qed{\hspace{12cm}$\diamond$}

   \def\diag{\mbox{\rm diag~}}

   \def\supp{\mbox{\rm supp}}
   \def\tr{\mbox{\rm tr~}}

   \def\({\langle}
   \def\){\rangle}
   \def\mb{\boldsymbol}

   \def\im{{\rm i}}

   \def\cA{{\mathcal A}}

   \def\cE{{\mathcal E}}
   \def\cF{{\mathcal F}}

   \def\cN{{\mathcal N}}

   \def\1{\mb1}

   \def\v0{{\bf 0}}

   \def\ov{\overline}

   \newtheorem{lemma}{Lemma}

   \newtheorem{example}{Ex.}
   \newtheorem{remark}{Remark}

\begin{document}

\title{Algebras of Interaction and Cooperation}
\author{Ulrich Faigle$^*$}\thanks{Mathematisches Institut, Universit\"at zu K\"oln, Weyertal 80, 50931 K\"oln, Germany\\ eMail: \texttt{faigle@zpr.uni-koeln.de}}

\date{\today}
\maketitle

{\bf Abstract:} Systems of cooperation and interaction are usually studied in the context of real or complex vector spaces. Additional insight, however, is gained when such systems are represented in vector spaces with multiplicative structures, \emph{i.e.}, in algebras. Algebras, on the other hand, are conveniently viewed as poly\-nomial algebras. In particular, basic interpretations of natural numbers yield natural polynomial algebras and offer a new unifying view on cooperation and inter\-action. For example, the concept of \emph{Galois transforms} and \emph{zero-dividends} of cooperative games is introduced as a nonlinear analogue of the classical Harsanyi dividends. Moreover, the polynomial model unifies various versions of Fourier transforms. Tensor products of polynomial spaces establish a unifying model with quantum theory and allow to study classical cooperative games as interaction activities in a quantum-theoretic context.

\medskip\noindent
{\bf Keywords:} {Activity system, cooperative game, decision system, evolution, Fourier transform, Galois transform, interaction system, measurement, quantum game.}

\section{Introduction}
What is a \emph{game}? While ''playing'' is often understood as an activity involving one or more human players, we take here the general game-theoretic approach of \cite{Faigle22} and think of an underlying \emph{system} that is characterized by the \emph{states} it can assume. The states typically depend on certain \emph{actions} that \emph{players} may or may not take in order to achieve other states. We do not necessarily require players to be human beings with human interests and feelings {\it etc.} and thus might refer to the players also simply as \emph{agents}. The agents my act, interact or cooperate according to  the rules within specified contexts.

\medskip
A mathematical analysis of games, of course, requires mathematical models for the underlying systems and their states. States have to be described as mathe\-matical objects. Similarly, observations on systems should be modeled accordingly as mathematical functions on the collection of states. The present paper con\-centrates on these aspects. In particular, we refer to \emph{quantum games} if the underlying systems fit or extend standard mathematical models of physical quantum systems.

\medskip
Questions about optimal strategies for the achievement of certain goals or about the existence of Nash equilibria {\it etc.} are disregarded. Given appropriate models, such questions lead to mathematical optimization problems that can be studied in their own right\footnote{See, \emph{e.g.}, \cite{Luenberger98,FaigleKernStill02,Faigle22}.}.

\medskip
In all of mathematical application analysis, and in game theory in particular, linear models have been proven to be of utmost importance. These models could be formulated and studied as abstract structures. The key in our analysis is the representation of relevant system parameters by polynomials as this point of view ties together many otherwise seemingly different models. Our polynomials are not \emph{polynomial functions} in \emph{variables} $x_i$ at the outset, but \emph{formal polynomials} in \emph{indeterminates} $x_i$. Appropriate interpretations of and substitutions into the in\-determinates then reveal fundamental aspects of various game-theoretic models.

\subsection{Outline of the presentation}
The algebraic structures in our analysis are suggested naturally when one represents system states as polynomials rather than vectors. The idea of expressing mathematical system characteristics {\it via} polynomials has a long tradition in general algebra (see, \emph{e.g.}, \cite{vdWaerden}). In fact, classical algebra arose from the wish to solve polynomial equations. The study of cooperative games in terms of an associated polynomial {\em function} in $n$ variables was initiated in \cite{Owen72}. Polynomial functions, however, typically imply just one particular algebraic structure: the addition and multiplication rules of scalar-valued functions.

\medskip
A substantially improved modeling flexibility is gained with \emph{formal} poly\-nomials rather than poly\-nomial functions. In this case, one deals with \emph{indeterminates} instead of \emph{variables}. Depending on the interpretation of the indeterminates, one is then led to various natural algebraic structures.

\medskip
The present approach is based on the three fundamental aspects of natural numbers: cardinalities of finite sets, binary representation of finite sets and representation of information in terms of $(0,1)$-sequences (Section~\ref{sec:natural-algebra}). Each of these aspects implies its own multiplication rule for formal polynomials.

\medskip
The polynomial model is then applied to cooperative games and activity systems. Section~\ref{sec:transforms} links polynomials to linear transforms and, in particular, to Fourier transforms. Moreover, the new concept of \emph{Galois transforms} is shown to appear naturally  and yields \emph{nonlinear} transforms of cooperative games. Interaction games are treated in Section~\ref{sec:interaction} with the emphasis on the role of a \emph{valuation} on a game-theoretic system as a particular state-dependent Heisenberg type measure on the system. A discussion with a perspective on future work concludes the presentation.

\subsection{Relation to earlier work}
Often, game theory is regarded as a scientific discipline in its own right\footnote{Also the fundamental mathematical theory~\cite{vNeumannMorgenstern44} proceeds this way.}. Moreover, cooperative games and strategic games with non-cooperative players  are treated separately. However, as questions about the computation of strategies became more and more of interest, the many connections of game theory with other mathematical fields ({\em e.g.}, mathematical optimization) became prominent. Computational questions have furthermore led to the emergence of game-theoretic research in computer science\footnote{\emph{cf.}~\cite{Nisan-etal07}.}. Moreover, a seemingly new area of game theory has been initiated where the games are supposed to be played according to the physical laws of quantum mechanics\footnote{See, \emph{e.g.,} \cite{PiotrowskiSladkowski03}}.

\medskip
As it turns out, a comprehensive approach to mathematical game theory is possible which ties together various areas of applied mathematics and includes the different aspects above as special cases\footnote{See \cite{Faigle22}.}. In particular, the relationship between game-theoretic cooperation and quantum mechanics has been recognized. In this sense, the present work is a continuation of the mathematical game theory research program begun in \cite{FaigleGrabisch16,FaigleGrabisch17}. As in the classical foundations of mathematical optimization\footnote{See, \emph{e.g.}, \cite{Luenberger98, FaigleKernStill02}.}, our mathematical model is essentially linear. However, it is observed that quadratic (and thus geometric, but nonlinear) measurement models arise naturally from the projection of linear operators onto lower dimensional spaces (\emph{cf.}  Section~\ref{sec:Heisenberg}).

\medskip
While  \cite{FaigleGrabisch24} explores some first aspects,  the present analysis pursues more generally the realization that the employment of polynomial algebras instead of pure abstract vector spaces as modeling tools offers a distinctly  more refined mathematical analysis and, furthermore, relates game theory to classical algebra.  No previous mathematical game-theoretical model seems to have taken this route before.

\section{Mathematical preliminaries}\label{sec:game-algebra}
Let $\N=\{0,1,2,\ldots,n,\ldots\}$ denote the set of natural numbers\footnote{{\em ''Die ganzen Zahlen hat der liebe Gott gemacht, alles andere ist Menschenwerk.''} (L. Kronecker) (\em ''The integers were created by the good Lord, everything else is the work of man.'') }. $\R$ is the field of real numbers and $\C$ the field of complex numbers, \emph{i.e.}, numbers of the form
$$
z = a +\im b \quad\mbox{where $a,b,\in \R$ and $\im^2 =-1$.}
$$
The complex number $z= a+\im b$ admits a representation in the form of de Moivre:
\begin{equation}\label{eq.de-Moivre}
	z = re^{i t}  = r(\cos t +\im \sin t) \quad\mbox{with real numbers $r,t\geq 0.$}
\end{equation}
The number $\ov{z} = a -\im b = r e^{-\im t}$ is the \emph{conjugate} of $z=a+\im b$ and has the property
$$
z\ov{z} = a^2 +b^2 = r^2  = |z|^2 .
$$

For arbitrary sets $X$ and $Y$, $X\times Y$ is the set of all pairs $(x,y)$ of elements $x\in X, y\in Y$. $Y^X$ is the collection of all functions $f:X\to Y$, which may be thought of as  \emph{valuations} of the elements of $X$ with values from $Y$.

The set $\C^X$ of all complex valuations of $X$ is a vector space under the usual addition and scalar multiplication of complex-valued functions.  The \emph{support} of $f \in \C^X$ is the set
$$
\supp(f) =\{x\in X| f(x)\neq 0\}.
$$
If $X=\{x_1,\ldots,x_n\}$ is a finite set, $\C^X$ may conveniently be identified with the $n$-dimensional coordinate space $\C^n$ and, similarly, $\C^{X\times X}$ with the space of all complex $n\times n$ matrices.

\subsection{Polynomials} A \emph{(formal) polynomial} relative to the set $X$ is an element $f\in \C^X$ with finite support. In order to emphasize the polynomial aspect, we write elements $f\in \C^X$ as formal sums
$$
      f = \sum_{x\in X} f_x x \quad\mbox{with the coefficients $f_x = f(x)$}.
$$
In the case $\supp(f)\subseteq  \{x_0,x_1,\ldots, x_n\}$,  we write the polynomial $f$ also in the form
$$
      f = \sum_{k=0}^n f_kx_k\quad\mbox{or, using superscripts $x^k \doteq x_k$,}\quad  f = \sum_{k=0}^n f_kx^k.
$$
We denote by $\C(X)$ the complex vector space of all polynomials in $\C^X$. The elements $x\in X$ are formally just \emph{indeterminates} without a numerical meaning in their own right. However, they can be given particular functional meaning in applications of the polynomial model.

\begin{remark} For notational convenience, we will often identify the indeterminate $x\in X$ with the polynomial $1 x\in \C(X)$ .
\end{remark}
\subsection{Polynomial functions}\label{sec:polynomial-functions}
If $X=\{x_0(t), x_1(t), x_2(t), \ldots\}$ is a family of complex-valued functions $x_k(t)$ in the variable $t$, one may think of a polynomial $p\in \C(X)$ as a complex-valued function in the variable $t$,
$$
   p(t) = \sum_{k\geq 0} c_k x_k(t),
$$
which arises as the corresponding linear combination of the functions in $X$. In the special case of the functions $x_k(t) = t^k$, $p(t)$ is a \emph{standard polynomial function}:
$$
   p(t) = c_0 +c_1x + c_2t^2 + \ldots + c_n t^n.
$$

\begin{remark}
The representation of polynomials as functions allows the application of the methods of differentiation and intergration in their analysis. For modeling purposes, however, it is important to retain the flexibility of formal polynomials.
\end{remark}

\begin{example}[Wave functions]\label{ex.wave-function} The choice $x_k(t) = (e^{2\pi \im t})^k = e^{2\pi\im kt}$ exhibits the \emph{wave functions} $v:\R\to \C$ of type
$$
v(t) = \sum_{k\geq 0} c_k[\cos(2\pi kt)] + \im\sin(2\pi kt)]
$$
as polynomial functions.
\end{example}

\subsection{Power series and generating functions}
In the case $X=\{x_0,x_1,x_2, \ldots,\}$ of indeterminates that are indexed by the natural numbers, we think of $f\in \C^X$ as a \emph{(formal power series)} with the notational representation
$$
        f = \sum_{k\geq 0} f_k x^k.
$$
Assuming $X=\{1, t,t^2,\ldots\}$ as a set of polynomial functions $x(t)=t^k$, on the other hand, the vector $f\in \C^X$ may define a complex function
$$
   f(t) = \sum_{k\geq 0} f_k t^k
$$
if the sum has a well-defined region of convergence in $\C$. In this case, the function $f(t)$ is interpreted as \emph{generating} the numerical parametes $f_k\in \C$ {\it via} evaluations of derivatives:
$$
    f_k = \frac{f^{(k)}(0)}{k!} \quad(k=0,1,2,\ldots).
$$
Generating functions have proven useful in probability theory or in the asymptotic analysis of combinatorial parameter sequences\footnote{See, \emph{e.g.}, \cite{Aigner97,Aigner07,Feller68}}. If $X$ is finite, a generating function is just a polynomial function.

\section{Incidence algebra} The \emph{function product} of two functions $f,g: X\to \C$ is the function $f\cdot g$ obtained from the usual product of the components
\begin{equation}\label{eq.function-product}
      (f\cdot g)(x) = f(x)g(x) \quad (x\in X) 
\end{equation}
and implies the corresponding multiplication not only of polynomials but of general power series. Note that the multiplication
$$
  (f,g) \mapsto f\cdot g
$$
is commutative and \emph{bilinear}, \emph{i.e.}, linear in each of the two components. All the multiplication structures we consider subsequently are derived from the function product, although not all are commutative.

\begin{remark}[Hermitian products]\label{rem.hermitian-products} Geometric considerations associate with a polynomial $f\in \C(X)$ a measure of its length as its \emph{norm}, namely the the real number
$$
   \|f\| = \sqrt{\sum_{x\in X} |f_x|^2} = \sqrt{\sum_{x\in X} f_x\cdot\ov{f_x}},
$$
which leads to the consideration of the \emph{hermitian} function product
$$
   (f,g) \mapsto f\cdot \ov{g},
$$
which is not linear but conjugated linear in its second component. The component sum of $f\cdot\ov{g}$ is the \emph{(hermitian) inner product}
\begin{equation}\label{eq.inner-product}
	(f|g) = \sum_x (f\cdot \ov{g})_x = \sum_x f_x\ov{g_x} \quad\mbox{with}\quad \|f\| = \sqrt{(f|f)}.
\end{equation}
\end{remark}

\subsection{Incidence functions}
For our puroses, an \emph{incidence function} on the set $X$ is a function 
$$
m:X\times X\times X\to \{0,1\}.
$$
 Given $m$, we
define the product $p = f\bullet g$ of the polynomials $f,g\in \C(X)$ as the polynomial with the coefficients
\begin{equation}
   p_x = \sum_{y,z\in X} m(x,y,z)  f_yg_z.
\end{equation}
The multiplication $(f,g)\mapsto f\bullet g$ on $\C(X)$ is no longer necessarily commutative, but still \emph{bilinear}.  We refer to $(\C(X), \bullet)$ as a \emph{(polynomial) incidence algebra}. This notion generalizes  Rota's~\cite{Rota64} fundamental algebraic model for the combinatorial analysis of discrete functions.
\begin{example}[M\"obius transform~\cite{Rota64}]\label{ex.Moebius-transform} Let $(X,\leq)$ be a partially ordered set and consider the incidence function
$$
   \zeta(x,y,z) = 1 \quad\Longleftrightarrow \quad x\geq y =z.
$$
Define the product $h=f\bullet g$ of the polynomials $f,g\in \C(X)$ as the polynomial with the coefficients

$$
    h_x = \sum_{y\leq x} f_yg_y = \sum_{y,z}\zeta(x,y,z)f_yg_z.
$$

The choice of unit coefficients $g_y = 1$, yields a linear operator
 $f\mapsto \zeta(f)$, where
$$
   \zeta(f)_x = \sum_{y\leq x} f_y.
$$
It is not difficult to see (and well-known) that the inverse operator $\mu=\zeta^{-1}$ exists. $w=\mu(f)$ is the \emph{M\"obius transform} of $f$ (relative to the given partial order on $X$) and has the property
$$
    f_x = \sum_{y\leq x} w_y \quad(x\in X).
$$
\end{example}

\begin{remark} In game theory, the coefficients $w_y$ of the M\"obius transform $w=\mu(f)$ are known as the \emph{Harsanyi~\cite{Harsanyi67} dividends} of the valuation $f$.
\end{remark}

If the incidence function $m$ is \emph{finitary}, \emph{i.e.}, if the set
$$
m (x)= \{(y,z)\in X\times X\mid m(x,y,z)=1\} \quad\mbox{is finite for all $x\in X$},
$$
the product rule of the algebra $(\C(X), \bullet)$, is easily extended to a product rule in the vector space $\C^*(X)$ of all linear functionals $L:\C(X)\to\C$. The product $F\bullet G$  is defined for  $F,G\in \C^*(X)$ as the linear functional  that acts on the basis polynomials $x\in X$ of $\C(X)$ as follows:
\begin{equation}\label{eq.linear-product}
(F\bullet G)(x) = \sum_{(y,z)\in m(x)} F_y G_z.
\end{equation}

\begin{remark}[Umbral calculus] Let $X=\{x^0,x^1,\ldots, \}$ be indexed by the natural numbers and consider the indicence function
	$$
	    m(x^k, x^i,x^j) =1  \quad \Longleftrightarrow\quad k= i +j.
	$$
\emph{Umbral calculus} is a technique to derive combinatorial identities by moving formally from superscripts to subscripts. It can be justified  {\it via} the application of linear operators and functionals on $\C(X)$ under the general multiplication rule
\begin{equation}\label{eq.linear-product}
	(f\bullet g)_k = \sum_{i\leq k}  \beta_{ik} f_j g_{k-i}.
\end{equation}
with suitable connection coefficients $\beta_{ik}$. See \cite{Rota-et-al73} for more details and many examples.
\end{remark}

\subsection{Monoidal algebras} The polynomial algebra $(\C(X),\bullet)$ is \emph{monoidal} if the multiplication turns $X$ into a \emph{monoid}, \emph{i.e}, if
\begin{enumerate}
\item $\1\bullet x = x\bullet\1 =x$ \;holds for some \emph{neutral element} $\1\in X$ and all $x\in X$;
\item $x\bullet (y\bullet z) = (x\bullet y)\bullet z \in X$ \;holds for all $x,y,z\in X$.
\end{enumerate}

Conversely, if $(X,\bullet)$ is a monoid, the multiplication extends to a multiplication of arbitrary polynomials in $\C(X)$ {\it via} the incidence function
$$
  m(x,y,z) =1 \quad\Longleftrightarrow\quad x= y\bullet z. 
$$

\begin{example}[Lattices]\label{ex.lattices} The partial order $L=(X,\leq)$ is a \emph{semilattice} if for all $x,y\in X$, there exists a unique largest element $x\wedge y$ such that
$$
   x\wedge y  \leq x \quad\mbox{and}\quad x\wedge y\leq y.
$$
If also a unique maximal element $x_{max}\in L$ exists, then $(X,\wedge)$ is a monoid with neutral element $x_{max}$. The semilattice $L$ is a \emph{lattice} if furthermore, any $x,y\in L$ admit a unique minimal element $x\vee y$ such that
$$
x\vee y\geq x \quad\mbox{and}\quad x\vee y\geq y,
$$
\emph{i.e}, if also the inverted partial order $L^*=(X,\geq)$ is a semilattice. If $L$ has a minimal element $x_{min}$, then $(X,\vee)$ is a monoid with neutral element $x_{min}$.
\end{example}

\section{Algebra of natural numbers}\label{sec:natural-algebra}
We denote the usual standard sum of natural numbers $i,j\in \N$ by $i+j$ and their product by $i\cdot j$ or simply $ij$. However, other important algebraic structures on $\N$ also offer themselves, depending on the representation of natural numbers.

\subsection{Binary representation and algebra}
Every natural number $ k < 2^n$ has a unique representation with $n$ binary digits $k_i$:
\begin{equation}\label{eq.binary-representation}
k = \sum_{i=0}^{n-1} k_i 2^k \quad(k_i\in \{0,1\}).
\end{equation}
So every $(0,1)$-string $\alpha \in \{0,1\}^\N$ with finite support  corresponds to a unique natural number
$$
   a = \sum_{i = 0}^\infty \alpha_i 2^i.
$$
Endowing the set $\{0,1\}$ with the binary addition rule
$$
   1\oplus 1 = 0,
$$
one obtains a binary addition rule for all natural numbers:
\begin{equation}\label{eq.binary-addition}
\big(\sum_{i = 0}^\infty \alpha_i 2^i\big)  \oplus \big(\sum_{i = 0}^\infty \beta_i 2^i\big) = \sum_{i = 0}^\infty (\alpha_i\oplus \beta_i) 2^i \quad(\alpha_i,\beta_i \in \{0,1\}),
\end{equation}
which is commutative and associative. Moreover, $0$ is the neutral element (as under the standard addition rule).

\medskip
Consider now a set $X=\{x^0,x^1,\ldots, x^n,\ldots \}$ of indeterminates $x^n$ with indices $n\in \N$. The binary addition suggests
a commutative and associative multiplication on $X$:
$$
    x^i \odot x^j = x^{i\oplus j}
$$
with the neutral element $\1=x^0$. The monoid $(X,\odot)$ implies the polynomial algebra $(\C(X),\odot)$ where the (binary) product of the polynomials $f,g\in \C(X)$ is the polynomial $f\odot g$ with the coefficients
$$
    (f\odot g)_n = \sum_{i\oplus j= n} f_i g_i.
$$

\medskip
\textsc{Impartial games.} We illustrate binary algebra with the example of the well-known $2$-person game \emph{Nim}, which involves two players and a finite set $N$ that is partitioned into $k$ blocks $N_i$ with $n_i=|N_i|$ objects each. A \emph{move} of a player consists in the choice of a non-empty block $N_i$ and the subsequent removal of at least one object from $N_i$. The players move alternatingly. A player loses if he cannot move on his turn. We associate with the Nim game its \emph{characteristic} polynomial as the binary monomial
$$
p(N)=x^{n_1}\odot x^{n_2}\odot \cdots\odot x^{n_k} = x^{n_1\oplus n_2\oplus \ldots +\oplus n_k}.
$$
By the Sprague-Grundy theory of combinatorial $2$-person games\footnote{See, \emph{e.g.}, \cite{Faigle22} for details} one then finds:
\begin{itemize}
\item The \emph{second} player has a winning strategy if $p(N)=\1\; (=x^0)$.
\item If $p(N) \neq \1$, the first player has a winning strategy.
\end{itemize}

\begin{remark} An actual winning strategy for the Nim game $N$ is easily computed: The current player moves, if possible, the game into a Nim situation $N'$ with characteristic polynomial $p(N') = \1$.
	
\medskip\noindent
Nim games are prototypical \emph{impartial} $2$-person games as each impartial game is strategically equivalent to a Nim game.
\end{remark}

\subsection{Set representations and Boolean algebra}\label{sec:boolean-algebra}
Each $(0,1)$-string $\alpha\in \{0,1\}^\N$ with components $\alpha_i$ describes a unique subset $A$ of $\N$ {\it via}
\begin{equation}\label{eq.indicator-representation}
  A = \supp(\alpha) = \{i\in \N\mid \alpha_i= 1\}.
\end{equation}
In fact, (\ref{eq.indicator-representation}) establishes a one-to-one correspondence between $\{0,1\}^\N$ and the subsets of $\N$. Moreover, strings with finite support correspond to finite subsets and, simultaneously, to natural numbers {\it via} their binary representation (\ref{eq.binary-representation}).

\medskip
Consider now an arbitrary $n$-element set $E=\{e_1,\ldots,e_n\}$. Each subset $A\subseteq E$ corresponds to a unique $(0,1)$-string $\alpha$ of length $n$,
$$
\alpha= \alpha_1\ldots\alpha_n \quad \mbox{with}\quad \alpha_i = 1 \quad\longleftrightarrow\quad e_i\in A,
$$
and also to a natural number
$$
 a = \sum_{i=1}^n \alpha_i 2^{i-1}  \;<\; 2^n .
$$
Consequently, the family $\cE$ of all subsets $A\subseteq E$ corresponds to the set
$$
\N_n =\{k\in \N\mid k< 2^n\}
$$
of natural numbers or to the family $\B_n =\{0,1\}^n$ of $(0,1)$-strings of length $n$. Hence, if $m=2^n$,  each polynomial of the form
$$
    p = p_0 x_0 +p_1 x_2 +\ldots + p_{m-1} x_{m-1}
$$
may equally well be indexed by the subsets of the $n$-element set $\cE$ or the binary $n$-strings:
$$
  p  =\sum_{k=0}^{m-1} p_kx_k \quad\longleftrightarrow \quad \sum_{A\in \cE} p_Ax_A \quad\longleftrightarrow \quad \sum_{\alpha\in \B_n} p_\alpha x_\alpha.
$$

The set-theoretic interpretation suggests a polynomial algebra based on Boolean lattice operations:
$$
    x_A \vee x_B = x_{A\cup B} \quad\mbox{and}\quad x_A\wedge x_B = x_{A\cap B}.
$$
$(X,\vee)$ is a monoid with neutral element $x_\emptyset$. $(X,\wedge)$ is a monoid with neutral element $x_E$ if  attention is restricted to polynomials indexed by the subsets of $E$.

\begin{example}[M\"obius algebra]\label{ex.Moebius-algebra} Let $X = \{x_A\mid A\subseteq E\}$ and set $q = \D\sum_{A\subseteq E} x_A$. Then 
$$
\sum_{T\subseteq S} p_T =(p\vee q)_S\quad\mbox{and}\quad \sum_{T\supseteq S} p_T =(p\wedge q)_S \quad\mbox{for any $p\in \C(X)$ and $S\subseteq E$}.
$$	
The example generalizes to lattices (see Ex.~\ref{ex.lattices}) in a straightforward manner.
\end{example}

\subsection{Concatenation and tensor products}\label{sec:tensor-product}
The \emph{concatenation} $\gamma=\alpha\star\beta$ of a string $\alpha\in \B_n$ with a string $\beta\in \B_m$ is the string
$$
\gamma= \gamma_1,\ldots \gamma_n\gamma_{n+1}\ldots \gamma_{n+m} \in \B_{n+m}
$$
where $\gamma_1\ldots \gamma_n = \alpha$ and $\gamma_{n+1}\ldots \gamma_{n+m} = \beta$. Its set-theoretic interpretation is straightforward:
\begin{itemize}
\item Let $N$ and $M$ be disjoint sets with $n=|N|$ and $m=|M|$ elements. Then $\gamma = \alpha\star \beta$ represents the set
$$
     A\cup B \subseteq N\cup M,
$$
if $\alpha$ represents the subset $A\subseteq N$ and $\beta$ the subset $B\subseteq M$.
\end{itemize}
Conversely, every string $\gamma\in \B_{n+m}$ decomposes into the product $\gamma=\alpha\star \beta$ of unique strings $\alpha\in \B_n$ and $\beta\in \B_m$.\footnote{However: if one does not refer explicitly to $\B_n$ and $\B_m$, the decomposition is \emph{not} unique.}

\medskip
The polynomial perspective on concatenation leads to the \emph{tensor product} of polynomials. Let $X$ and $Y$ be two sets of indeterminates and denote by $X\otimes Y$ the collection of all pairs $x\otimes y$ of indeterminates $x\in X$ and $y\in Y$.

\medskip
We define the \emph{tensor product} of $f\in \C(X)$ and $g\in \C(Y)$, in analogy with (\ref{eq.function-product}), as the formal product
\begin{equation}\label{eq.tensor-product}
	\left(\sum_{x\in X} c_x x \right)\;\otimes\; \left(\sum_{y\in Y} d_y y\right) \;=\;  \sum_{(x,y)\in X\times Y} c_xd_y (x\otimes y).
\end{equation}
The tensor product is bilinear but, in contrast to the function product,  not necessarily commutative. The collection of all tensor products is not closed under taking sums and hence not a vector space in its own right. The tensor products generate the vector space
$$
\C(X)\otimes \C(Y) = \C(X\otimes Y) \;\; (\;\simeq\; \C(X\times Y)\,).
$$

\subsection{Entanglement}
If $f\in \C(X)$ and $g\in \C(Y)$ are polynomials of norm $\|f\|=1=\|g\|$, then the polynomial $f\otimes g\in \C(X\otimes Y)$ has norm $\|f\otimes g\|=1$. Intuitively, $f\otimes g$ arises as the \emph{independent} join of $f$ and $g$. Notice however:
\begin{itemize}
\item There may be polynomials $h\in \C(X\otimes Y)$ that cannot be \emph{disentangled}, \emph{i.e.},
$$
   h\neq f\otimes g \quad\mbox{holds for all $f\in \C(X)$ and $g\in \C(Y)$}.
$$
\end{itemize}

\medskip
This fact has an important consequence in mathematical modeling. In models where systems states are represented by polynomials, states of the tensor product may possibly not be split into independent components.

\begin{remark} States that cannot be disentangled may appear counterintuitive in physical system models and have received considerable attention particularly in quantum theory\footnote{See, \emph{e.g}, \cite{Bell64}. }.
\end{remark}

\section{Cooperative games}
The model of a {\em cooperative TU-game} $\Gamma=(N,v)$ involves a set $N=\{e_1,\ldots,e_n\}$ of $n$ \emph{players}. A subset $S\subseteq N$ is a \emph{coalition}. Letting $\cN$ be the family of all coalitions (subsets) of $N$, a valuation $v:\cN\to\R$ is assumed to be given with $v(S)$ as the value generated by the coalition $S$. We define the \emph{characteristic} polynomial of $\Gamma$ in its set-theoretic notation with indeterminates $x_S$ as
$$
p^v= \sum_{S\in \cN} v_S x_S \quad\mbox{where \quad $v_S = v(S).$}
$$
For its numerical interpretation, recall that the members of the set $\cN$ correspond to the $m=2^n$ natural numbers $0\leq s< m$. So the characteristic polynomial takes the form
$$
   p^v \quad\longleftrightarrow\quad \sum_{s=0}^{m-1} v_s x^s
$$
 with the index numbers $s$ corresponding to the coalitions $S\in \cN$.

\medskip
Finally, denoting by $|\sigma\)$ the indeterminate $x_S$ corresponding to the coalition $S\in \cN$ in its representation as a $(0,1)$-string $\sigma\in \B_n$, we have
$$
  p^v \quad\longleftrightarrow \quad\sum_{\sigma\in \B_n} v_\sigma|\sigma\).
$$

\begin{remark} Interpreting the valuation $v$ of the cooperative game $\Gamma=(N,v)$ as a function on $\B_n=\{0,1\}^n$, $v$ is also known as a \emph{pseudo-boolean function}\footnote{See, \emph{e.g.}, \cite{BorosHammer02}.}.  Moreover, $v$ is usually called the \emph{characteristic function} of $\Gamma$. In the present context, however, it is useful to allow the latter terminology with a greater flexibility (see the next Section~\ref{sec:characteristic-functions}).
\end{remark}

\begin{remark} A more general model of \emph{cooperative games with restricted cooperation}\footnote{See \cite{Faigle89}.} considers only valuations $v:\cF\to\R$ on a given subfamily $\cF\subseteq\cN$ of coalitions. The ''augmenting systems'' of \cite{Bilbao03} are special cases of this model, for example. We leave it to the interested reader to extend the present theory accordingly.
\end{remark}

\subsection{Characteristic functions}\label{sec:characteristic-functions}
There are various functions associated with cooperative games (relative to the set $N$ of players) which result from appropriate substitutions into the characteristic polynomial $p^v$. Examples from Section~\ref{sec:polynomial-functions} are the polynomial function $p^v(t)$ and the wave function $v(t)$, both with one real variable $t$.

\medskip
If a complex variable $t_S$ is substituted for each indeterminate $x_S$, the characteristic polynomial $p^v$ yields a linear functional $\hat{v}: \C^\cN\to \C$ with values
$$
\hat{v}(t) =\sum_{S\in \cN} v_S t_S  \quad(t\in \C^\cN).
$$

\begin{example}[Multinomial extensions] Owen's~\cite{Owen72}\label{ex.Owen}   \emph{multinomial extension}  of the valuation $v:\cN\to \R$ is obtained if $n$ real variables $t_1,\ldots, t_n$ are associated with the $n$ players $e_1,\ldots,e_n$ and the polynomial functions
$$
   t_S = \prod_{e_i\in S}t_i \prod_{e_j\notin S}(1-t_j)
$$
are substituted into the indeterminates $x_S$ of the characteristic polynomial:
$$
p^v(t_1,\ldots, t_n) = \sum_{S\in \cN} v_S   \prod_{e_i\in S}t_i \prod_{e_j\notin S}(1-t_j).
$$	

\end{example}

Similar polynomials are obtained if one considers the set $\C$ of complex numbers (or the set $\R$ of real numbers) to be endowed with a lattice structure $L = (\C,\vee,\wedge)$ and replaces the usual scalar addition and multiplication by the lattice operations:
$$
     c + d \mapsto c\vee d\quad\mbox{and}\quad c\cdot d \mapsto c\wedge d.
$$

\begin{example}[Sugeno integrals] Assume $L=(\R,\wedge,\vee)$ to be a lattice and define for a valuation $v:\cN\to \R$ its  \emph{lattice polynomial} as
$$
p^v = \bigvee_{S\in \cN} v_S\wedge x_S.
$$
Associating real variables $t_i$ with the elements $e_i$ of $N$ as in Ex.~\ref{ex.Owen}, the substitution of
$$
   t_S = \bigwedge_{e_i\in S} t_i
 $$	
then yields the lattice polynomial function
$$
p^v(t) = \bigvee_{S\in \cN} v_S\wedge t_S \quad(t\in \R^\cN).
$$
The \emph{Sugeno integral}\footnote{See \cite{Grabisch2016,Marichal09} for details.} relative to a valuation $v:\cN\to \R$ is the lattice polynomial with the lattice operations
$$
     c \vee d = \max\{c,d\}  \quad\mbox{and}\quad c\wedge d= \min\{c,d\}.
$$
\end{example}

\subsection{Activity systems}\label{sec:activity-systems}
Suppose that the players $i\in N$ in the cooperative game $\Gamma=(N,v)$ decide to become \emph{active} in $\Gamma$ with respective probabilities $\pi_i$. Then the evaluation
$$
p^v(\pi_1,\ldots, \pi_n) = \sum_{S\in \cN} v_S \pi_S \quad\mbox{with $\pi_S = \D\prod_{i\in S}\pi_S$}
$$	
of the multinomial extension yields the expected value of the ensuing cooperation -- provided the players take their individual decisions \emph{independently} from each other\footnote{The \emph{cooperative fuzzy games} of \cite{Aubin81} make a similar independence assumption on the players.}. In general, \emph{i.e.}, without the independence assumption, the expected value results from the substitution of some probability distribution $\pi$ on the coalition family $\cN$ into the indeterminates of the characteristic polynomial:
\begin{equation}\label{eq.expected-value2}
	E^v(\pi) = \sum_{S\in \cN} v(S)\pi(S).
\end{equation}
So $\pi$ describes probabilistically a \emph{state} the cooperative $\Gamma$ is in and the value it is expected to generate in this state.

\medskip
As a probabilistic augmentation of the cooperation model, let us think of the members of  $N$ simply as \emph{agents} (players, or coalitions of players, or physical entities {\it etc.}) and define an \emph{activity system}\footnote{According to the standard interpretation of quantum mechanics,  finite-dimensional quantum systems are activity systems in the present sense. (See also Section~\ref{sec:Heisenberg}).} $\cA_N$ relative to a finite set $N=\{1,\ldots,n\}$ of agents with collection $\cN$ of subsets as a system whose states are described by complex vectors $u\in \C^\cN$ of unit length $\|u\|=1$ or, equivalenty, (squared) norm
$$
\|u\|^2 = \sum_{S\in \cN} |u_S|^2 = 1.
$$

The (real) parameters $|u_S|^2$ specify a probability distribution on $\cN$. Hence, if $v:\cN\to\R$ is a measuring function on $\cN$, its expected observation value will be
\begin{equation}\label{eq.expected-value2}
	E^v(u) = \sum_{S\in \cN} v(S)|u_S|^2
\end{equation}
if $\cA_N$ is in the activity state $u$. Conversely, every probability distribution $\pi$ on $\cN$ arises from some unit length vector $u^\pi$. For example:
$$
u^\pi = (\sqrt{\pi_S}: S\in \cN) \in \C^\cN.
$$

\begin{remark}[Hermitian multiplication]\label{rem.Hermitian-multiplication} The states of $\cA_N$ correspond to poly\-nomials $s \in \C(X)$ with norm $\|s\| =1$. The hermitian product $(s,s) =s\cdot \ov{s}$ has the components $(s,s)_T = |s_T|^2$. So the expected value results from the evaluation of the substitution of $v$ into $(s,s)$.

\end{remark}

\section{Transforms}\label{sec:transforms}
A \emph{transform} on $\C(X)$ is an operator $T:\C(X)\to \C(X)$. For example, the M\"obius transform (Example~\ref{ex.Moebius-transform}) $\mu$ is a linear transform. In game theory, the underlying partial order is usually assumed to be set-theoretic:
$$
    x_A \leq x_B  \quad\Longleftrightarrow\quad A\subseteq B,
$$
and yields a Haranyi-type valuation $w=\mu(v)\in \R^\cN$ for each $v:\cN\to \R$ such that
$$
     v(B) =\sum_{A\subseteq S} w(A)\quad\mbox{for all $B\in \cN$.}
$$

The idea of the characteristic polynomial suggests also structurally important algebraic transforms on $\C(X)$ that are nonlinear and have not received attention in game theory so far.

\medskip
Throughout this section, we assume an underlying family $\cN$ of $m=2^n$ coalitions $S$ of an $n$-element player set $N$, which are
represented by the natural numbers $0\leq s < m$.

\subsection{Fourier transforms} For a given scalar $\omega\in \C$, the basic $\omega$-transform of an indeterminate $x$ is its scalar adjustment
$$
    x \mapsto \omega x.
$$

\subsection{The discrete Fourier transform}
In the case of the standard polynomial algebra $(\C(X),\cdot)$, the indeterminates have the product decomposition
$$
    x^s = x^1\cdot x^1 \cdots x^1 \quad\mbox{($s$ times)}.
$$
So $x^1 $ is the basic indeterminate. The $\omega$-transform $x^1 \mapsto \omega x^1$ implies
$$
  (\omega x^1)^s = \omega^s x^s \quad\mbox{for all $s$}
$$
and is represented by its characteristic polynomial
$$
   f^{(\omega)} =\sum_{k = 0}^{m-1} \omega^k x^k.
$$
The associated \emph{$\omega$-transform} $F^\omega$ on the whole space $\C(X)$ is the linear transform
$$
     p \;\mapsto\; F^\omega(p)=f^\omega\cdot p.
$$
Since $F^\omega$ is linear, it suffices to keep track of its effect on the basic polynomials $x^s$:
$$
    F^\omega(x^s) = f^\omega \cdot x^s =  \sum_{k=0}^{m-1} \omega^{sk} x^k \quad\mbox{for all $s=0,1,\ldots,m-1$}.
$$
For the choice $\omega = e^{2\pi \im/m}$, $F^\omega$ is known as the \emph{discrete Fourier transform} on $\C(X)$ and yields for $v\in \C^\cN$,
$$
F^\omega(v_0x^0+  v_1 x +\ldots+v_{m-1}x^{m-1}) = \hat{v}_0x^0 + \hat{v}_1 x +\ldots +\hat{v}_{m-1} x^{m-1}
$$
with the \emph{Fourier coefficients}
$$
     \hat{v}_j = \sum_{k=0}^{m-1} v_k e^{2\pi\im jk/m}  \quad(j=0,1,\ldots,m-1).
$$	

It is easily verified that the discrete Fourier transform has a Fourier type inverse:
$$
{v}_k = \frac1m\sum_{j=0}^{m-1} \hat{v}_j e^{-2\pi\im jk/m}  \quad(k=0,1,\ldots,m-1).
$$

\begin{remark} Viewed alternatively, the Fourier coefficients $\hat{v}_j$ are obtained by either substituting the interpolation parameters $t_j=e^{2\pi\im j/m}$ into the corresponding polynomial function or by the substitution of the parameters $j$ into the corresponding wave function.
\end{remark}

\subsection{The quantum Fourier transform} The so-called \emph{quantum Fourier transform} $QFT_n$ is the discrete Fourier transform in disguise when one replaces the indeterminates $x^s$  by the corresponding binary bit strings $|s\)\in \B_n$. The effect on the string $|s\) \in \B_n$ is:
$$
QFT_n|s\) = \sum_{k=0}^{m-1} \omega^{sk} |k\) =  \sum_{k=0}^{m-1} e^{2\pi \im sk/m}|k\).
$$

\medskip
\begin{remark} The quantum Fourier transform is a major tool in the theory of quantum computation\footnote{see, \emph{e.g.}, \cite{Nielsen-Chuang00} }.
\end{remark}

\subsection{ The parity transform} In the polynomial algebra $(\C(X),\odot)$ with the binary multiplication
$$
    x^s\odot x^k = x^{s\oplus k} \quad\mbox{for all $s,k\in \B_n$}\; ,
$$
the $n$ indeterminates $x_1,\ldots,x_n$ associated with the $n$ unit elements in $\B_n$ are basic in the sense
$$
x_s = \prod_{i:s_i=1} x_i \quad\mbox{and hence}\quad \prod_{i:s_i=1}(\omega x_i) = \omega^{|s|} x_s,
$$
where $|s| = s_1 + \ldots +s_n$ is the usual sum of the components of the vector $s\in \{0,1\}^n$. In set-theoretic notation, the analogue of the polynomial $f^\omega$ is the polynomial
$$
   h^\omega = \sum_{S\in \cN} \omega^{|S|} x_S
$$
which implies the transform
$$
     p \;\mapsto h^\omega\odot p.
$$
The choice $\omega = -1$ yields the \emph{parity transform}
$$
    p\;\mapsto\; \hat{p} = h^{-1}\odot p.
$$
with the \emph{binary Fourier coefficients}
$$
     \hat{v}_S = \sum_{T\subseteq N}(-1)^{|S\cap T|}v_T \quad\mbox{for any $v=\D\sum_{S\subseteq N}  v_S x_S\in \C(X)$.}
$$

\begin{remark} The parity transform is also called \emph{binary Fourier transform}\footnote{For interesting applications of the binary transform in game theory, see, \emph{e.g.}, \cite{Kalai02,Odo}. }. The \emph{Hadamard transform} $H$ is the normalized parity transform:
$$
    H(v) = \frac{1}{\sqrt{m}}\hat{v} \quad\mbox{with the self-inverse property}\quad H(H(v)) =v.
$$
\end{remark}

\subsection{Galois transforms and zero-dividends} Recall from the Fundamental Theorem of Algebra\footnote{See, \emph{e.g.,} \cite{vdWaerden}} that the polynomial function
$$
v(t) = v_0 + v_1 t +v_2 t^2 +\ldots +v_{k-1}t^{k-1} + t^k
$$
with $k\geq 1$, admits a parameter vector $g=(g_1,\ldots,g_k)\in \C^k$ such that
\begin{equation}\label{eq.product-representation}
 v(t) = (t+g_1)(t+g_2) \cdots (t + g_k).
\end{equation}
One thus has
\begin{equation}\label{eq.symmetric-functions}
	\begin{array}{cclll}
	v_{k-1}&=& g_1+g_2+\ldots+g_k\\
	v_{k-2} &=& g_1g_2+\ldots +g_2g_3+ \ldots + g_{k-1}g_nk\\
	v_{k-3} &=& g_1g_2g_3 + g_1g_2g_4+\ldots\\
	&\vdots && \\
	v_0 &=& g_1g_2\cdots g_k .
\end{array}
\end{equation}

\begin{remark} The numbers $z_i=-g_i$ are exactly the zeroes of the polynomial function $v(t)$. So the parameter vector $(g_1,\ldots,g_k)$ is unique up to permutations of coordinates.
\end{remark}

If $v(t)$ is the characteristic polynomial function of a valuation $v\in \C^m$, we augment $g$ to a vector in $\C^m$:
$$
  g(v) = (g_1,\ldots,g_k,\ldots,g_m) \quad\mbox{with $g_\ell = 1$ for $\ell \geq k$}.
$$
We call the valuation $g(v)$ a \emph{Galois transform} of $v$.

\medskip
If $\Gamma=(N,v)$ is a cooperative game with characteristic polynomial function
$$
v(t) = v_0 + v_1 t +v_2 t^2 +\ldots + v_kt^k
$$
and $v_k\neq 0$ we \emph{normalize} $v(t)$ to the polynomial $\tilde{v}(t) = v(t)/v_k$ with the leading coefficient $\tilde{v}_k = 1$ and call the parameters $z_i = -g_i$ relative to the Galois transform $g(\tilde{v})$ of the normalized game $\tilde{\Gamma}=(N,\tilde{v})$ the \emph{zero-dividends} of $\Gamma$.

\medskip
Notice that the normalized game $\tilde{\Gamma}$ can be reconstructed from its Galois transform.

\begin{example}[Zero-normalized games]
	Cooperative game theory often makes the assumption that a cooperative game $(N,v)$ is \emph{zero-normalized}, \emph{i.e.}, has the property $v(\emptyset) = 0$. Since $\tilde{v}(\emptyset)$ is the product of the zero-dividends $g_i$ of $v$, one has
	$$
	v(\emptyset) = 0  \quad\Longleftrightarrow\quad g_1g_2\cdots g_{m}=0
	$$
	Hence $v$ is zero-normalized if and only if $g_i = 0$ holds for some zero-dividend $g_i$ of $v$.
	
\end{example}

\section{Interaction games}\label{sec:interaction}
Consider a set $N=\{1,\ldots, n\}$ of general entities $i$ and specify an \emph{interaction game} on $N$ as a valuation
$$
V:N\times N\to \R,
$$
where $V_{ij} = V(i,j)$ is the \emph{interaction worth}\footnote{See, for example, the {\em interaction indices} in various application modeling contexts~\cite{Beliakov19,sciencedirect}.} of $i,j\in N$. Where $X$ is a set of $n$ indeterminates $x_i$, the associated characteristic polynomial is
$$
\chi^V = \sum_{i=1}^n\sum_{j=1}^n V_{ij} (x_i\otimes x_j),
$$
which means that interaction games refer to the tensor product space $\C(X\otimes X)$ (and not to the vector space $\C(X)$) or, more precisely, to the real tensor space $\R(X\otimes X)$. Abstractly, $V$ is given as a real $n\times n$ matrix with the coefficients $V_{ij}$.

\medskip
Let $A\in \R^{N\times N}$ be an arbitrary real matrix whose coefficients $A_{ij}$ have the interpretation of enhancing the interaction worth relative to $V$:
\begin{itemize}
	\item If $i$ and $j$ interact at level $A_{ij}$, their interaction produces the value $V_{ij}A_{ij}$.
\end{itemize}
Hence the interaction instance $A$ produces the game's overall value as
\begin{equation}\label{eq.inner-product0}
	\chi^V(A)=  \sum_{i,j\in N} V_{ij} A_{ij} = \tr(A^TV)
\end{equation}
where the \emph{trace} $\tr(C)$ of a matrix $C$ is the sum of its diagonal coefficients. In other words, $\chi^V(A)$ equals the usual euclidean inner product $(V|A)$  of the two matrices $V,A$, considered as $n^2$-dimensional parameter vectors. This fact
yields the dual interpretations:

\begin{enumerate}
	\item {\em An interaction game $V$ is a linear functional on the space of all interaction instances $A$.}
	\item {\em An interaction instance $A$ is a linear functional on the space of all interaction games $V$.}
\end{enumerate}

\subsection{The hermitian perspective}
Setting $A^+ = (A+A^T)/2$, one finds that a matrix $A\in \R^{n\times n}$ decomposes into a symmetric matrix $A^+$ and a skew-symmetric matrix $A^-$:
\begin{equation}\label{eq.symmetry}
	A = A^+ + A^- \quad\mbox{where}\quad (A^+)^T = A^+, (A^-)^T = - A^-.
\end{equation}
Moreover, is it straightforward to check that the symmetry decomposition (\ref{eq.symmetry}) of $A$ is unique. Associate now with $A\in \R^{n\times n}$ the well-defined \emph{hermitian} matrix
$$
\hat{A} = A^+ +\im A^- \in \C^{n\times n}
$$
and let $\H_n\subseteq \C^{n\times n}$ be the family of all hermitian $n\times n$ matrices.

\medskip
Clearly, $\H_n$ is isomorphic to $\R^{n\times n}$ with respect to the field $\R$ of real scalars. ($\H_n$ is not a complex vectors space, however.) Recall that the \emph{adjoint} $C^*$ of the complex matrix $C\in \C^{n\times n}$ is the transpose $\ov{C}^T$ of the conjugated matrix $\ov{C}$ and note

\medskip
\begin{lemma}\label{l.hermitian} For any matrices $C\in \C^{n\times n}$ and $A,B\in \R^{n\times n}$, one has
	\begin{enumerate}
		\item $C\in \H_n$ if and only if $C^* = C$ (\emph{i.e.}, $C$ is self-adjoint).
		\item $(A|B) = \tr(B^TA) = \tr(\hat{B}\hat{A})=(\hat{A}|\hat{B})$.
	\end{enumerate}
\end{lemma}

\noindent
Hence:

\begin{itemize}
	\item {\em Interaction games can equally well be studied within the context of the hermitian matrix space $\H_n$.}
\end{itemize}

\medskip
The adjoint $u^*$ of the (column) vector $u\in \C^n$ is the row vector of conjugated components of $u$. So $u u^*$ is a self-adjoint $n\times n$ matrix. Indeed, an important  general characterization of self-adjoint matrices is:

\medskip
\begin{lemma} The matrix $C\in \C^{n\times n}$ is self-adjoint if and only if there are real scalars $\lambda_k$ and vectors $u_k\in \C^n$ such that
	\begin{equation}\label{eq.hermitian}
		C = \sum_k \lambda_k u_k u_k^*.
	\end{equation}
\end{lemma}

\Pf Because the matrices $u_ku^*_k$ are self-adjoint, a matrix $C$ of the form (\ref{eq.hermitian}) is self-adjoint. To see the converse, recall that $\R^{n\times n}$ and $\H_n$ are isomorphic. So it suffices to consider $(0,1)$-matrices $A\in \R^{n\times n}$ with exactly one non-zero entry $A_{ij} = 1$. Moreover, we may assume $n=2$.

\medskip
If $A$ is diagonal, one has $\hat{A}=A$, for which the claim is obviously true. So assume, for example,
$$
A = \begin{pmatrix} 0 &1\\ 0 &0 \end{pmatrix}\quad\mbox{and hence}\quad C = \hat{A} = \frac12\begin{pmatrix} 0 & 1+\im\\  1-\im &0\end{pmatrix}
$$
with the two real eigenvalues $\lambda_1=+\sqrt{2}$ and $\lambda_2 = -\sqrt{2}$. Consequently,  the $2$-dimensional space $\C^2$ admits a basis of eigenvectors of $C$, which imply the claim.

\qed

\begin{remark} The spectral decomposition of self-adjoint matrices\footnote{See, \emph{e.g.}, \cite{Weidmann12}.} shows that the vectors $u_k$ in (\ref{eq.hermitian}) can be chosen as eigenvectors of $C$ with real eigenvalues $\lambda_k$ also in the case $n>2$.
\end{remark}

\medskip\noindent
Consider, for example, an activity system $\cA_N$ as in Section~\ref{sec:activity-systems} relative to the valuation $v:\cN\to \R$. Let $u\in \C^\cN$ be a state vector of norm $\|u\|=1$ with the associated self-adjoint $m\times m$ matrix $U = u u^*$ of complex coefficients
$$
U_{ST} = u_S\ov{u_{T}}.
$$
Where $V=\diag(v)$ is the diagonal interaction matrix with diagonal coefficients $V_{SS} = v(S)$, one now observes:
$$
(V|U) = \sum_{S\subseteq N} v(S)\ov{u_S}u_S =  \sum_{S\subseteq N} v(S)|u_S|^2=  E^v(u).
$$

\medskip
So $\cA_N$ is recognized as the restriction of an interaction game on $\cN$ with characteristic function $V=\diag(v)$ to interaction instances of the form $u u^*$.

\subsection{Heisenberg measurements}\label{sec:Heisenberg}
A \emph{measurement operator} on the interaction system $N$ is a functional
$$
\mu:\R^{N\times N} \to \R
$$
which produces the measurement result $\mu(A)$ if the system is in a state that corresponds to the activity instance $A$.  If the functional $\mu$ is linear, the measurement acutally represents an interaction game. So there exists a matrix $M\in \R^{n\times n}$ such that
$$
   \mu(A) = (M|A) =(\hat{M}|\hat{A}),
$$
which means that $\mu$ can be understood as a linear functional on the (real) vector space $\H_n$ of all hermitian matrices. Assuming
$$
   \hat{A} = \sum_k \lambda_k u_ku_k^*
$$
for suitable real parameters $\lambda_k$ and (complex) vectors $u_k\in \C^n$, one has
\begin{equation}\label{eq.Heisenberg-measurement}
   \mu(A) = \sum_k \lambda_k \mu(u_ku_k^*) =\sum_k \lambda_k (\hat{M}|u_ku_k^*).
\end{equation}

A \emph{Heisenberg measurement operator}\footnote{''Heisenberg'' refers to the fact that this measurement model is standard in quantum theory.}  on the (complex) vector space $\C^n$ is a real-valued functional  $\gamma:\C^n \to\R$ of the form
$$
\gamma(u) = u^*Gu = (G|uu^*)\quad\mbox{with $G\in \C^{n\times n}$}.
$$
Since $\gamma$ is real-valued, one may assume that $G$ is self-adjoint. Notice that a Heisenberg measurement operator is not linear. On the other hand, (\ref{eq.Heisenberg-measurement}) shows:
\begin{itemize}
\item {\em A Heisenberg measurement operator on $\C^n$ arises from the restriction of a \emph{linear} measurement operator on the interaction system $N$ to interaction instances of the hermitian form $uu^*$.}
\end{itemize}
Similarly, the linear measurement matrix $\hat{M}$ of the operator $\mu$ in (\ref{eq.Heisenberg-measurement}) is of the form
$$
\hat{M} = \sum_j \delta_j w_jw_j^*
$$
for suitable real numbers $\delta_j$ and vectors $w_j\in \C^n$, which implies
\begin{equation}\label{eq.interaction-measurement}
	\mu(A) = \sum_i\sum_j \lambda_i\mu_j (w_jw_j^*|u_iu_i^*).
\end{equation}

Consequently, the fundamental linear interaction measures are recognized to be of the form
\begin{equation}\label{eq.inner-product}
	\mu_w(u) = (ww^*|uu^*) = |w^*u|^2 \quad\mbox{with}\quad w,u\in \C^n.
\end{equation}

\section{Discussion}
Linearity plays an important role in mathematical application models. The mathe\-matical analysis, however, will reveal more characteristic features of the underlying systems when the model is not just considered to be a vector space with scalar multiplication, but an algebra, \emph{i.e}, additionally equipped with a multiplication operation for vectors. Suitable multiplication operations are naturally associated with multiplication rules for polynomials, which renders polynomial models powerful and flexible.

\medskip
The present approach shows that mathematical models for cooperation and interaction connect with important aspects of classical algebra and combinatorics. For example, the representation of coalitions by natural numbers embeds the representation of cooperative games into the context of Galois theory, \emph{i.e.}, the theory of solving algebraic equations. Future work of exploring this area of mathematical system analysis in more detail appears to be promising.

\medskip
The polynomial model also underlines the aspect of quantum-theoretic models as interaction systems and, conversely, embeds cooperation and interaction into the setting of physical quantum systems. The evolution of such systems can be mathematically understood in a far broader context (see, {\em e.g.}\cite{FaigleGierz18}).
Moreover, the apparatus of theoretical physics can be brought to bear on general systems of cooperation and interaction. In particular, Hamiltonians of cooperative games can be expected to provide considerable insight into fundamental laws according to which such systems behave.

\end{document}